\documentclass[12pt,a4paper]{amsart}
\usepackage[T1]{fontenc}
\usepackage{lmodern}
\usepackage{amsmath,amssymb,amsthm}
\usepackage{microtype}
\usepackage[a4paper,margin=18mm]{geometry}
\usepackage[hidelinks]{hyperref}
\newcommand{\R}{\mathbb R}
\newcommand{\Tr}{\operatorname{Tr}}
\newcommand{\dd}{\,\mathrm d}
\newtheorem{theorem}{Theorem}

\author{Fabio Nicola}
\address{Dipartimento di Scienze Matematiche ``G. L. Lagrange'',
Politecnico di Torino, Torino, Italy}
\email{fabio.nicola@polito.it}

\title{On the concentration of the ambiguity function on centered balls}
\begin{document}
\begin{abstract} We show that the Gaussian uniquely maximizes the $L^2$ concentration of the ambiguity function on every centered ball, up to time-frequency shifts and multiplication by a unimodular constant. The proof relies on ideas from quantum information theory.
\end{abstract}
\maketitle

\section{Statement and proof of the main result} 
We follow the conventions of \cite{G} concerning the short-time Fourier transform and the Wigner distribution: for $f,g\in L^2(\R^d)$, $x,\omega\in\R^d$,
\[
 V_gf(x,\omega)=\int_{\R^d}f(t)\overline{g(t-x)}e^{-2\pi i\omega\cdot t}\dd t,
 \qquad
 W(f)(x,\omega)=\int_{\R^d}f(x+t/2)\overline{f(x-t/2)}e^{-2\pi i\omega\cdot t}\dd t.
\]
Let $B_R\subset\R^{2d}$ be the ball of radius $R$ centered at the origin and let
$\varphi(t)=2^{d/4}e^{-\pi|t|^2}$, with $t\in\mathbb{R}^{d}$, be the standard normalized Gaussian. 

Lieb and Ostrover~\cite{LO} proved that, for each
$R>0$, the functional 
\[
\int_{B_R}W(f)(z)\dd z
\]
 is maximized by $\varphi$ among all
$f\in L^2(\R^d)$ with $\|f\|_2=1$ (here $z=(x,\omega)$, and $dz$ stands for the Lebesgue measure in $\R^{2d}$). 

Gr\"ochenig and Neuhauser~\cite{GN} asked whether a similar result holds for the functional 
\[
\int_{B_R}|V_ff(z)|^2\dd z,
\]
which measures the \(L^2\) concentration of the ambiguity function, or equivalently, the mass of the spectrogram of \(f\) on \(B_R\). This question has been circulating in the time-frequency analysis community for several years. We show that the answer is affirmative.

\begin{theorem}\label{mainteo} Let $R>0$. For every $f\in L^2(\R^d)$, with $\|f\|_{L^2}=1$, we have 
\[
 \int_{B_R}|V_f f(z)|^2\dd z\leq \int_{B_R}e^{-\pi|z|^2}\dd z.
 \]
 
Equality holds if and only if $f(t)=e^{i\theta}e^{2\pi i t \cdot \omega_0} \varphi(t-x_0)$ for some
$\theta\in\R$ and $x_0,\omega_0\in\R^d$. 
\end{theorem}
\begin{proof}
The desired conclusion will follow from the following facts. 

\begin{itemize}
\item[(a)] The Lieb--Ostrover bound~\cite{LO} mentioned above extends to every density operator
$\rho=\sum_{j\ge0}p_jP_{\psi_j}$, with $p_j\ge0$, $\sum_jp_j=1$ and $\|\psi_j\|_{L^2}=1$, where $P_\psi f=\langle f,\psi\rangle\psi$. Indeed,
$W(\rho)(z)=\sum_jp_jW(\psi_j)(z)$, with convergence in $L^2(\R^{2d})$, in particular in $L^1(B_R)$.

\item[(b)] Writing $\check f(t)=f(-t)$, the covariance of the Wigner distribution and the Moyal identity
give (see, e.g., \cite[(2) and (9)]{Janssen1982}) 
\[
 |V_ff(z)|^2=\bigl(W(f)*W(\check f)\bigr)(z),\qquad z\in\R^{2d}.
\]

\item[(c)] For density operators $\rho_1,\rho_2$ in $L^2(\R^d)$, their \emph{symmetric quantum convolution} $\rho_1\boxplus \rho_2$
is a density operator satisfying, in our normalization,
\[
 \rho_1\boxplus \rho_2=\Tr_2\!\left[U_{1/2}(\rho_1\otimes \rho_2)U_{1/2}^{*}\right],
 \qquad
 W(\rho_1\boxplus \rho_2)(z)=2^d\bigl(W(\rho_1)*W(\rho_2)\bigr)(\sqrt2\,z);
\]
see \cite[(37), (47)]{BDLR}, where $U_{1/2}$ is the so-called
beam-splitter unitary. In fact, we only need {\em the existence} of a density operator whose Wigner distribution is $2^d\bigl(W(\rho_1)*W(\rho_2)\bigr)(\sqrt2\,z)$.  
\end{itemize}

Thus, with $\rho_f=P_f\boxplus P_{\check f}$,
$W(\rho_f)(z)=2^d|V_ff(\sqrt2\,z)|^2$, and point~(a) yields
\[
 \int_{B_R}|V_ff(z)|^2\dd z
 =\int_{B_{R/\sqrt2}}W(\rho_f)(z)\dd z
 \le\int_{B_{R/\sqrt2}}W(\varphi)(z)\dd z
 =\int_{B_R}e^{-\pi|z|^2}\dd z.
\]
It remains to characterize the cases of equality. The uniqueness assertion in
\cite{LO}, applied to the spectral decomposition of $\rho_f$, forces
$\rho_f=P_\varphi$; hence $|V_ff(z)|^2=e^{-\pi|z|^2}$. The marginal identities for the spectrogram then imply that equality holds in Heisenberg's uncertainty inequality in each coordinate, which yields the asserted form of $f$. Conversely, every such function gives equality by the invariance of \(|V_ff|^2\) under time-frequency shifts.
\end{proof}

\textit{Remark.}  Alternatively, for the point~(c) one could check that the function \[
F_\sigma W(\rho_1)(z/\sqrt2)F_\sigma W(\rho_2)(z/\sqrt2)\]
 satisfies the conditions of the quantum Bochner Theorem (see e.g.~\cite{dG,Holevo}), using the Schur product theorem ($F_\sigma$ denoting the symplectic Fourier transform). 

The analogous Gaussian maximization result fails for general powers \(|V_ff|^p\). For instance, it fails for \(p=1\), as can be seen by considering the first Hermite function.

\section*{Acknowledgments} 
I would like to thank Karlheinz Gröchenig for suggesting the problem addressed in this note, for sharing the manuscript \cite{GN}, and for his comments on a preliminary version of this note.

\end{document}